\documentclass[12pt]{article}
\usepackage{setspace}
\usepackage[hmargin=3 cm,vmargin={2.5 cm,2.5 cm}]{geometry}
\usepackage{amsmath}
\usepackage{apacite}
\usepackage[authoryear]{natbib}
\usepackage{amssymb}
\usepackage{graphicx,multirow}
\usepackage{algorithm}
\usepackage{algorithmic}
\usepackage{subfigure}
\usepackage{multicol,float}
\usepackage{bbm}
\usepackage{url}
\usepackage{tikz,pgfplots}
\usepackage{amsthm}
\usepackage{authblk}
\usepackage[pdfborder={0 0 0}]{hyperref}
\newcommand{\ind}{\mathbbm{1}}

\newcommand{\Keywords}[1]{\par\noindent{\small{\em Keywords\/}: #1}}

\def\EE{\mathbb{E}}
\def\PP{\mathbb{P}}

\def\L{{\cal L}}

\begin{document}
\onehalfspacing

\title{Managing losses in exotic horse race wagering}

\author[1]{Antoine Deza}
\author[2]{Kai Huang}
\author[3]{Michael R. Metel}

\affil[1]{Advanced Optimization Laboratory, Department of Computing and Software, McMaster University, Hamilton, Ontario, Canada\\
\url{deza@mcmaster.ca}}
\affil[2]{DeGroote School of Business, McMaster University, Hamilton, Ontario, Canada\\
\url{khuang@mcmaster.ca}}
\affil[3]{Laboratoire de Recherche en Informatique, Universit\'e Paris-Sud, Orsay, France\\
\url{metel@lri.fr}}

\maketitle

\begin{abstract}
We consider a specialized form of risk management for betting opportunities with low payout frequency, presented in particular for exotic horse race wagering. An optimization problem is developed which limits losing streaks with high probability to the given time horizon of a gambler, which is formulated as a globally solvable mixed integer non-linear program. A case study is conducted using one season of historical horse racing data.
\end{abstract}


\section*{Introduction}

Since the mid 1980's, horse racing has witnessed the rise of betting syndicates akin to hedge funds profiting from statistical techniques similar to high frequency traders in the stock market \citep{kaplan2002}. This is possible as parimutuel wagering is employed at racetracks, where money is pooled for each bet type, the racetrack takes a percentage, and the remainder is disbursed to the winners in proportion to the amount wagered.\\

Research on horse racing stems in large part due to the fact that it can be viewed as a simplified financial market. Research on important economic concepts such as utility theory \citep{weitzman}, the efficient market hypothesis \citep{asch}, and rational choice theory \citep{rosett} can be done in a straight forward manner, given horse racing's discrete nature, fixed short term contract lengths and attainable sets of historical data for empirical study.\\

Optimization in the horse racing literature can be traced back to \citet{isaacs1953} deriving a closed form solution for the optimal win bets when maximizing expected profit. \citet{hau81} utilized an optimization framework to show inefficiencies in the place and show betting pools using win bet odds to estimate race outcomes. In particular, they used the \citet{kelly1956} criterion, maximizing the expected log utility of wealth, and found profitability when limiting betting to opportunities where the expected return was greater than a fixed percentage. More recently, \citet{smoc2010} derived a simple procedure for optimal win bets under the Kelly criterion through analysis of the Karush-Kuhn-Tucker optimality conditions.\\

Having found a favourable opportunity in a gambling setting, such as betting on the outcome of flipping a biased coin, the Kelly criterion answers the question of how much to wager. For example, if the probability of heads is $\PP(H)=0.6$, with even payout odds, and wealth $w$, we can determine how much to wager on heads, $x$, by maximizing the expectation of the log of our wealth after the toss, $\max\limits_{x}0.6\log(w+x)+0.4\log(w-x)$, which has an optimal solution of $x^*=0.2w$, telling us to always wager 20\% of our current wealth. Kelly style betting is widely recognized both in academia~\citep{maclean2011} and in practice, being used professionally in blackjack~\citep{carlson2001}, general sports betting~\citep{wong2009}, and in particular horse race betting~\citep{wong11}. Positive aspects of the Kelly criterion are that it asymptotically maximizes the rate of return of one's wealth, and assuming one can wager any fraction of money, it never risks ruin. The volatility of wealth through time is too large for most though, as $\PP(w_t\leq \frac{w_0}{n}|t>0)\approx \frac{1}{n}$~\citep{thorp2006}, e.g. there's approximately a 10\% chance your wealth in the future will be 10\% of what it currently is using the Kelly criterion. As a result, many professional investors choose to employ a fractional Kelly criterion \citep{thorp2008}, which has been shown to possess favourable risk-return properties by \citet{maclean1992}, with betting half the Kelly amount being popular amongst gamblers~\citep{Pound05}.\\

There are several different types of wagers one can place on horses, including what are known as exotic wagers, which include the exactor, triactor and superfecta, which require the bettor to pick the first two, three and four finishers in order, respectively. The exotic wagers are popular among professional gamblers, as superior knowledge of the outcome of a race is better rewarded, and the more exotic the bet, the higher the advantage one can attain~\citep{ben08}. For this reason we focus on the superfecta bet, the most exotic wager placed on a single race.

\section*{Time horizon}
\label{sec:OM}

In recognition of the similarities between parimutuel horse race betting and financial markets, we see superfecta betting being most similar to the purchase of deep out of the money options, with the general trend of a successful strategy being small steady losses through time with infrequent large gains. Speaking of his experience as a key member of a Hong Kong horse racing gambling syndicate, \citet{wong11} states that investing in horse racing is more stressful than in the stock market, and that for professional groups wagering in exotic pools it is normal not to have a winning wager once in three months. Once the losing streak terminates a large profit is achieved, but in the interim, there will be various sources of pressure. Doubt in the system may set in leading to the potential for irrational decisions to be made, based not on statistical findings but emotion.\\

It would be ideal to have a mechanism to control losing streaks,
not only to avoid failure but to determine if a losing streak is in range with the current strategy or if an investigation into the system is warranted. As this is a form of risk management, we consider such methods from stock portfolio management. The most famous framework is mean-variance portfolio optimization based on the work of ~\citet{mark}, where one maximizes the expected return subject to a constraint which limits the variance in portfolio returns. One of the criticisms of this model is that the use of variance as a measure of risk penalizes both positive and negative deviations in the same manner. Given the expected positive skewness of superfecta returns this would be particularly problematic for our application.\\

A popular risk measure proposed to replace variance is the value at risk (VaR) \citep{brandi}, which estimates the maximum amount a portfolio could lose over a given time period at a given confidence level $1-\alpha$. Maximizing the Kelly criterion subject to a VaR constraint has been considered previously by \citet{maclean2004} in the context of allocating investment capital to stocks, bonds and cash over time. Let $S$ represent the set of top four horse finishers with each $s\in S$ corresponding to a sequence of 4 horses, with $x=\{x_s\}$ being our decision variables dictating how much to wager on each outcome $s$, and $P(x)$ being the random payout given our decision vector $x$. Let the outcome probability of $s$ be denoted as $\pi_s$, with $\pi_{x}=\sum_{s\in S}\pi_s\ind_{\{x_s>0\}}$ being the probability of having a winning bet.
We can now limit our betting strategy's VaR to be no greater than $v$ by enforcing the chance constraint $\PP(P(x)-\sum_{s\in S}x_s\geq -v)\geq 1-\alpha$. More broadly, chance constrained optimization enables the accommodation of data uncertainty by enforcing affected constraints with a given probability. For more background, see \citet{shapiro2009}.\\

VaR calculations typically use a small $\alpha$, being concerned with large potential losses near the tail of the distribution. Tail risk is not a concern in our setting as the most that could possibly be lost is the amount we wager, which we expect to occur most of the time, in fact, a VaR constraint with $v>0$ in our setting corresponds to a betting limit for $\alpha<1-\pi_x$.\\

Though risk measures concerning tail losses seem unapplicable, a VaR constraint with $v\leq0$ enables the control of losing streaks. Let $\tau$ be the gambler's time horizon, for which we desire to set as the limit for potential losing streaks with high probability.
For a betting decision $x$, let $\tilde{\pi}_{x}=\PP(P(x)-\sum_{s\in S}x_s\geq -v)$ and $B_{x}\sim \text{binomial}(\tau,\tilde{\pi}_{x})$ be the random number of times at least $-v$ dollars is earned repeating the race $\tau$ times with the same wager $x$. In order to enforce the gambler's time horizon, we require that $\PP(B_{x}\geq 1)\geq 1-\alpha$, which implies $\tilde{\pi}_{x}\geq 1-\alpha^{\frac{1}{\tau}}$. Assuming independence between races, limiting betting decisions to those which have a VaR $\leq 0$ with confidence of at least $1-\alpha^{\frac{1}{\tau}}$ ensures that a non-negative return on a race will occur with a probability of at least $1-\alpha$ over the next $\tau$ races.

\section*{Optimization model}

A conceptual optimization model is displayed below.
Using the Kelly criterion, the objective is to maximize the expected log of wealth, where $w$ is the current wealth of the gambler. We also simplify our VaR constraint notation, enforcing from now on $\PP(P(x)-\sum_{s\in S}x_s\geq v)\geq 1-\alpha^{\frac{1}{\tau}}$ for $v\geq 0$.

\begin{alignat}{6}
&\max&&\text{ }\EE\log(P(x)+w-\sum_{s\in S}x_s)\nonumber\\
&\mbox{s.t. }&&\sum_{s\in S}x_s\leq w \nonumber\\
&&&\PP(P(x)-\sum_{s\in S}x_s\geq v)\geq 1-\alpha^{\frac{1}{\tau}}\nonumber\\
&&&x_{s}\geq 0\hspace{5 pt} s\in S \nonumber
\end{alignat}

\section*{Case study}
\label{sec:CS}

The optimization model was tested using historical race data from the 2013-2014 season at Flamboro Downs, Hamilton, Ontario, Canada. This amounted to a total of 1,168 races. Race results, including the payouts, pool sizes, and final win bet odds were collected from \citet{TI}. Handicapping data, generated by \citet{CB}, was collected from \citet{HPI}. The first $70\%$ of the race dataset was used to calibrate the race outcome probabilities and payout model, with the remaining $30\%$ of races used for out of sample testing.

\subsection*{Estimating outcome probabilities and payouts}
\label{sec:RO}

The multinomial logistic model, first proposed by \citet{Bolt86}, is the most widely used method of estimating the probability of each horse winning a race. Given a vector of handicapping data on each horse $h$, $v_{h}$, the horses are given a value $V_{h}=\beta^Tv_{h}$, and assigned winning probabilities $\pi_{h}=\frac{e^{V_{h}}}{\sum_{i=1}^{n}e^{V_i}}$. A three factor model was used, including the log of the public's implied win probabilities from the win bet odds, $\log{\pi^p_h}$, and the log of two CompuBet factors. The analysis was performed using the {\it mlogit} package~\citep{croi12} in {\it R}. Details of the handicapping data and the statistical estimation can be found in the subsection {\it Estimating win probabilities} in the appendix.\\

The ~\cite{Harv73} model assumes the probability that a horse finishes $m^{th}$ equals the probability that it wins against the horses that didn't finish $1^{st},..,{m-1}^{th}$.
The conditional probabilities are $\pi_{ij|i}=\frac{\pi_j}{1-\pi_i}$, $\pi_{ijk|ij}=\frac{\pi_k}{1-\pi_i-\pi_j}$, and $\pi_{ijkl|ijk}=\frac{\pi_l}{1-\pi_i-\pi_j-\pi_k}$, where for example, $\pi_{ijk|ij}$ is the probability estimate of horses $i$, $j$, and $k$ finishing first, second, third, given horses $i$ and $j$ finished first and second. Multiplying together with $\pi_i$, $\pi_{ijkl}=\frac{\pi_i\pi_j\pi_k\pi_l}{(1-\pi_i)(1-\pi_i-\pi_j)(1-\pi_i-\pi_j-\pi_k)}$. This model was found to be biased towards favourite horses by ~\cite{Lo94b} and ~\cite{Lo94}.
We use the improved approximation derived by ~\cite{Lo08},
$\pi_{ijkl}=\pi_i\frac{\pi^{\lambda_1}_j}{\sum_{s\neq i}\pi^{\lambda_1}_s}\frac{\pi^{\lambda_2}_k}{\sum_{s\neq i,j}\pi^{\lambda_2}_s}\frac{\pi^{\lambda_3}_l}{\sum_{s\neq i,j,k}\pi^{\lambda_3}_s}$, where $\lambda_1$, $\lambda_2$ and $\lambda_3$ are calibrated to the historical race data. As the log-likelihood is separable, optimal $\lambda_i$'s were determined individually using multinomial logistic regression. The results of the statistical estimation can be found in the subsection {\it Estimating superfecta probabilities} in the appendix.\\

The superfecta payout function for sequence $s$ is approximately $P_{s}(x)=x_{s}\frac{(Q+\sum_{u\in S}x_u)(1-t)}{Q_{s}+x_{s}}$, where $Q$ is the superfecta pool size, $Q_{s}$ is the amount wagered on sequence $s$ by other gamblers, and $t$ is the track take. The payout per dollar wagered is typically rounded down to the nearest nickel, termed breakage, but this is unlikely to be significant and is omitted from the formula. The only information available to bettors is the value of $Q$. The approach taken to estimate $Q_s$ is motivated by the work of ~\citet{kant08} who fit the win probabilities of the Harville model to the money wagered on quinella bets using multinomial maximum likelihood estimation. Let $b$ be the minimum allowable bet, with larger wagers being a multiple of $b$. The amount wagered on sequence $s$ is $Q_{s}=\frac{Q(1-t)}{P_s}$, where $P_{s}$ is the amount paid on a \$1 wager. Let $n=\frac{Q_s}{b}$ be the number of bets placed on $s$ out of $N=\frac{Q}{b}$, which we assume follows a binomial distribution.
We model the public's estimate of superfecta outcome probabilities using a discount model with the public's implied win probabilities, so for $s=\{i,j,k,l\}$,
$\pi^p_s=\frac{(\pi^p_i)^{\theta_1}}{\sum_h(\pi^p_h)^{\theta_1}}\frac{(\pi^p_j)^{\theta_2}}{\sum_{h\neq i}(\pi^p_h)^{\theta_2}}\frac{(\pi^p_k)^{\theta_3}}{\sum_{h\neq i,j}(\pi^p_h)^{\theta_3}}
\frac{(\pi^p_l)^{\theta_4}}{\sum_{h\neq i,j,k}(\pi^p_h)^{\theta_4}}$.
Let $\pi^p_{s,u}$ and $\pi^p_{s,l}$ represent the numerator and denominator of $\pi_s^p$.
The likelihood function, using data from $R$ historical races assumed to be independent, with $w_r$ being the winning sequence in race $r$, is
$\L(\theta)\propto\Pi_{r=1}^R(\pi^p_{w_r})^{n_r}(1-\pi^p_{w_r})^{N_r-n_r}$. The log-likelihood is a difference of concave functions,
$\log\L(\theta)\propto\sum_{r=1}^Rn_r\log(\pi^p_{w_r,u})+(N_r-n_r)\log(\pi^p_{w_r,l}-\pi^p_{w_r,u})-N_r\log(\pi^p_{w_r,l})$. This function was minimized twice using {\it fminunc} in Matlab, the first with an initial guess that the public uses the Harville model, $\theta_i=1$, the second assuming that the public believes superfecta outcomes are purely random, $\theta_i=0$, with both resulting in the same optimal solution. Statistical estimation results can be found in the subsection {\it Estimating public's superfecta probabilities} in the appendix.\\

In our simulations we use a point estimate of $Q_s$. Ideally, we want to be highly certain that our profit will be at least equal to $v$ for scenarios satisfying the chance constraint. In an attempt to achieve this, we approximate the uncertainty of our estimate of the public's superfecta probabilities by modeling $\theta\sim N(\hat{\theta},\Sigma)$, where $\hat{\theta}$ is our maximum likelihood estimate and $\Sigma$ is our estimated covariance matrix of $\theta$, taken as the inverse of the observed Fisher information. For each race wagered on, we took 9999 samples of $\theta$ and generated a sample of $\pi^p$ from each. We then took $\pi_s^p$ as the $99^{th}$ sample percentile by setting it to its $9900^{th}$ ordered statistic, and set $Q_s=\pi_s^pQ$.

\subsection*{Optimization formulation}

We now formulate the optimization program as it will be solved, assuming we want to use a fractional Kelly strategy with fraction $f$. This is accomplished by multiplying the optimal solution by $f$, then rounding each bet to the closest multiple of $b$ to generate a valid wager. The chance constraint is implemented using binary variables $z_s$, which indicate that a bet will be placed on outcome $s$, which should generate a profit of at least $v$. This can be modeled as $b\lfloor \frac{f}{b}x_{s}\rceil\frac{(Q+\sum_{u\in S}b\lfloor \frac{f}{b}x_{u}\rceil)(1-t)}{Q_{s}+b\lfloor \frac{f}{b}x_{s}\rceil}-\sum_{u\in S}b\lfloor \frac{f}{b}x_{u}\rceil \geq (v+w)z_s-w$. In order to preserve convexity, we use the approximation $x_{s}\frac{(Q+\sum_{u\in S}x_{u})(1-t)}{Q_{s}+x_{s}}-\sum_{u\in S}x_{u}\geq (\frac{v}{f}+w)z_s-w$, which ignores the rounding and the non-linearity of the payoff function. We also require that if $z_s=1$, then $x_s\geq \frac{b}{f}$ to ensure a wager will be placed on outcome $s$, which can be implemented by the constraint $\frac{b}{f}z_s\leq x_s$. We then enforce the chance constraint by $\sum_{s\in S}\pi_sz_s\geq 1-\alpha^{\frac{1}{\tau}}$. Note that our implementation is an approximation of $\tilde{\pi}_{x}\geq 1-\alpha^{\frac{1}{\tau}}$ given our point estimate of $Q_s$.

\begin{alignat}{6}
\label{eq:F1}
&\max&&\text{ }\sum_{s\in S}\pi_{s}\log(x_{s}\frac{(Q+\sum_{u\in S}x_u)(1-t)}{Q_s+x_{s}}+w-\sum_{u\in S}x_u)\\
&\mbox{s.t. }&&\sum_{s\in S}x_s\leq w \nonumber\\
&&&\sum_{s\in S}\pi_sz_s\geq 1-\alpha^{\frac{1}{\tau}} \nonumber\\
&&&x_{s}\frac{(Q+\sum_{u\in S}x_u)(1-t)}{Q_{s}+x_{s}}-\sum_{u\in S}x_u \geq (\frac{v}{f}+w)z_s-w\hspace{5 pt}s\in S \nonumber\\
&&&\left(\frac{Q_s+\frac{b}{f}}{Q_s}\right)^{z_s}\leq\frac{Q_s+x_s}{Q_s}\hspace{11 pt}s\in S \nonumber\\
&&&z_{s}\in \{0,1\}\hspace{82 pt}s\in S \nonumber\\
&&&x_{s}\geq 0\hspace{103 pt}s\in S \nonumber
\end{alignat}

The objective function of (\ref{eq:F1}) is not concave and the third constraint is not convex. We use the 1 to 1 mapping proposed by \citet{Kall08}, $y_s=\log(x_s+Q_s)$, resulting in the following program which is convex after relaxing the binary constraints on $z_s$. We have written the constraints $\frac{b}{f}z_s\leq x_s$ equivalently above as $\left(\frac{Q_s+\frac{b}{f}}{Q_s}\right)^{z_s}\leq\frac{Q_s+x_s}{Q_s}$ in order to achieve convex constraints after the change of variable.

\begin{alignat}{6}
\label{eq:F2}
&\max&&\text{ }\sum_{s\in S}\pi_{s}\log(Q+w-(t+(1-t)Q_se^{-y_s})\sum_u e^{y_u})\\
&\mbox{s.t. }&&\sum_{s\in S}e^{y_s}\leq w+Q \nonumber\\
&&&\sum_{s\in S}\pi_sz_s\geq (1-\alpha^{\frac{1}{\tau}}) \nonumber\\
&&&Q-(t+(1-t)Q_se^{-y_s})\sum_u e^{y_u}\geq (\frac{v}{f}+w)z_s-w\hspace{5 pt}\forall s\in S \nonumber\\
&&&z_s\ln\left(\frac{Q_s+\frac{b}{f}}{Q_s}\right)\leq y_s-\log{Q_s}\hspace{11 pt}s\in S \nonumber\\
&&&z_{s}\in \{0,1\}\hspace{110 pt}s\in S \nonumber\\
&&&y_s\geq \log(Q_s)\hspace{99 pt}s\in S \nonumber
\end{alignat}

\subsection*{Implementation}

All computation was conducted on a Windows 7 Home Premium 64-bit, Intel Core i5-2320 3GHz processor with 8 GB of RAM, in Matlab R2016a using OPTI toolbox v2.16. For each race, IPOPT~\citep{IPOPT06} was first used to solve (\ref{eq:F2}) without the time horizon constraint. If $\sum_{s\in S}x_s=0$, we do not bet on the current race and if $\sum_{s\in S}\pi_sz_s\geq 1-\alpha^{\frac{1}{\tau}}$ we take the result as the solution. If $\sum_{s\in S}x_s>0$ but $\sum_{s\in S}\pi_sz_s < 1-\alpha^{\frac{1}{\tau}}$, we proceed to solve the full problem using Bonmin's~\citep{Bonmin2008} {\it B-Hyb} algorithm. None of the default stopping criteria was altered in OPTI's optimization settings, so the maximum execution time was limited to 1,000 seconds, the maximum number of iterations to 1,500 and the maximum function evaluations to 10,000. With these settings it was not always guaranteed that the optimal solution was found. In order to improve solution quality, we only considered a subset of possible outcomes to wager on.
Outcomes were ordered by probability times profit from placing a single wager of $\frac{b}{f}$ on each,
$\pi_{s}\frac{b}{f}\left(\frac{(Q+\frac{b}{f})(1-t)}{Q_s+\frac{b}{f}}-1\right)$, with the top 50\% of outcomes considered in the optimization program. In our dataset, the estimated probability of these outcomes had a median value of $84\%$ and always contained the winning outcome.

\subsection*{Results}
\label{sec:R}

\begin{figure}[H]
 \centerline{
\resizebox{1.0\textwidth}{!}{
\begin{tikzpicture}
\scriptsize
    \begin{axis}[xlabel=Race,ylabel=Wealth,
                legend style={legend pos=outer north east}]
                \addplot[mark=none,mark size=0.75,draw=black, thick]
                table[x=x,y=y1]
            {results3500.dat};
            \addplot[mark=none,mark size=0.75,draw=darkgray, thick]
                table[x=x,y=y2]
            {results3500.dat};
            \addplot[mark=none,mark size=0.75,draw=gray, thick]
                table[x=x,y=y3]
            {results3500.dat};
            \addplot[mark=none,mark size=0.75,draw=lightgray, thick]
                table[x=x,y=y4]
            {results3500.dat};
    \legend{$\tau=\infty$,\hspace{5.5mm}$40$,\hspace{5.5mm}$30$,\hspace{5.5mm}$20$}
    \end{axis}
\end{tikzpicture}}}
\caption{Wealth over the course of 350 races at Flamboro Downs with $v=0$.} \label{T9}
\end{figure}
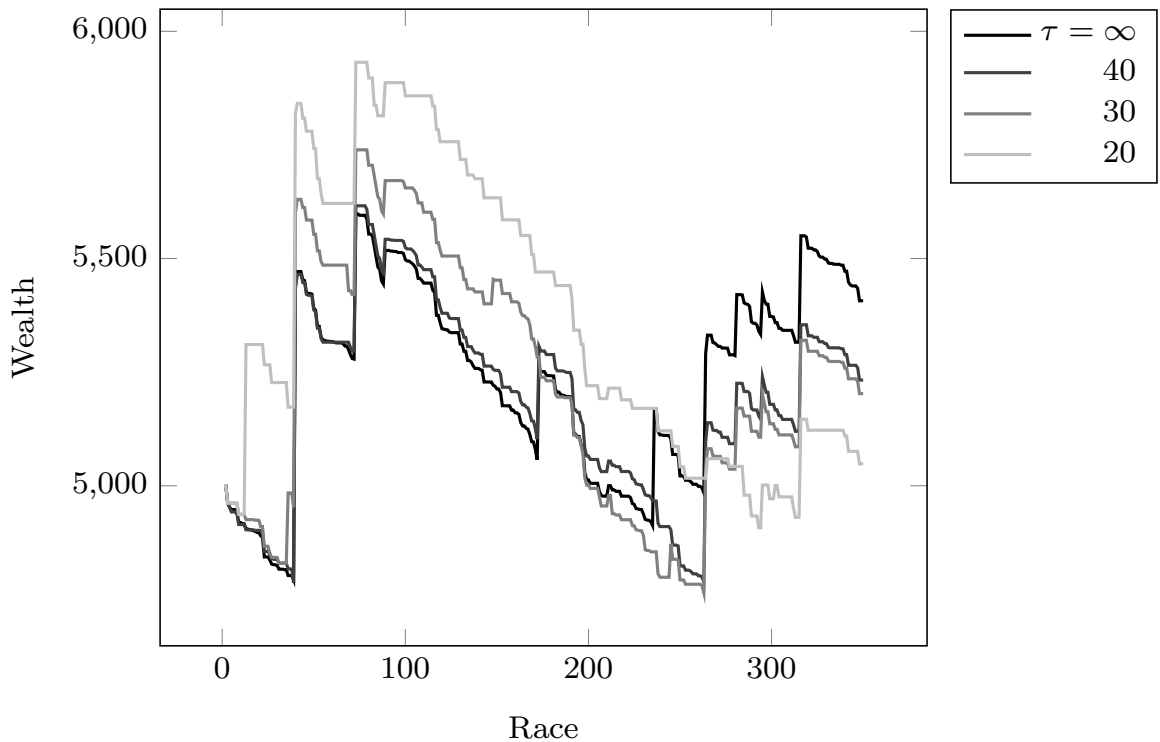

Testing was done on a total of 350 races from Flamboro Downs, where $t=24.7\%$ and $b=\$0.2$. Due to its success in practice, we set $f=0.5$.  Given our optimal betting solution, the realized payout was calculated by adjusting the published payout to account for our wagers and breakage. Four simulations were done with the gambler's initial wealth set to \$5,000. The wealth through time for all are plotted in Figure \ref{T9}, with statistics displayed in Table \ref{T5}. A preliminary simulation was done with $\tau=\infty$. The longest losing streak was found to be $52$ races. Given this number, simulations were done with $\tau=40$, $30$ and $20$, with $\alpha=0.05$. We set $v=\$20.24$, which was the minimum positive profit achieved in a race with $\tau=\infty$, and the maximum value of $v$ for which the longest losing streak with $\tau=\infty$ remains unchanged, while also ensuring that negligible winning bets do not end losing streaks for other values of $\tau$. In our simulations, we considered a losing streak to end after a profit $P(x)-\sum_{s\in S}x_s>0.99v$ was realized.\\

 \begin{table}[H]
\centerline{
\resizebox{0.75\textwidth}{!}{
\renewcommand{\arraystretch}{1}
\begin{tabular}{lccccccc}
   $\tau$&Loss streak&Total return (\%)&Races bet&Bet per race\\
  \hline
$\infty$&$52$&$7.8$&$224$&$10.8$\\
40 &$40$&$4.5$&$163$&$14.4$\\
30 &$27$&$4.0$&$123$&$20.5$\\
20 &$19$&$1.0$&$61$&$37.9$\\
\hline
\end{tabular}}}
\caption{Optimization results} \label{T5}
\end{table}

Examining Table \ref{T5}, Loss streak is the maximum losing streak over races bet on, Total return is the total return over the 350 races, Races bet is the total number of races bet on, and Bet per race is the average bet per race. The length of losing streaks were successfully limited to the chosen time horizon, but we can see there is a trade off between risk and return, resulting in a reduction in profit using the chance constrained model. The chance constraint forced us to be more selective in which races we wagered on, and increased the average amount bet per race as it became required to be profitable in more outcomes.

\section*{Conclusion and future research}
\label{sec:C}
We have developed a methodology for limiting losing streaks given a gambler's time horizon through the use of chance constrained optimization, exemplified in exotic horse race wagering. Initial results using one season of historical racing data have been presented which show the viability of the method by effectively limiting losing streaks for different chosen time horizons. Certain approximations were used which could be addressed in future research. Point estimates of outcome probabilities, $\pi_s$, as well as the amount wagered on each outcome by the public, $Q_s$, were utilized. Taking into further account the uncertainty of these estimates could improve results. Though the focus of this work has been on horse racing, we feel this general methodology could be applicable to any gambling or investing setting which have low probability outcomes with high payouts, such as investing in deep out of the money options.

\section*{Acknowledgements}

The authors thank the anonymous referees for their valuable comments. This work was partially supported by the Natural Sciences and Engineering Research Council of Canada Discovery Grant programs (RGPIN-2015-06163, RGPIN06524-15), and by the Digiteo Chair C\&O program.

\bibliographystyle{apacite}
\bibliography{CCHorseRacingReferences}

\section*{Appendix}

\subsection*{Estimating win probabilities}
{
A number of factors and their logarithms were considered, displayed in Table \ref{T2} below. The domain of each factor is listed in brackets, but all were  normalized to be between 0 and 1 for statistical use. The first six factors are from \citet{CB}, with the other two from the race program and result.
\begin{table}[H]
\centerline{
\begin{tabular}{p{2.4cm}p{0cm}p{11.6cm}}
  \hline
	\multicolumn{3}{c}{\bf{Factor description}} \\
  \hline
  Post && Starting position of the horse (1-9).\\
  Pre && The quality of the data available for each horse (30-100).\\
  Form && The overall success of this horse in recent starts (10-130).\\
  Class && The horse's performance relative to the class of its competition in recent races (52.8-95).\\
  Speed && An adjusted speed rating using the daily track variant, track condition, and the track-to-track speed variant (113.3-128.1 seconds).\\
  Driver Points && The driver's rating (4-39).\\
  $\pi^{ML}_{h}$ && The winning probability implied by the morning line odds.\\
  $\pi^m_{h}$ && The winning probability implied by the final winning bet odds.\\
\end{tabular}}
\caption{Win probability considered factors.} \label{T2}
\end{table}
}
Systematically removing the least significant factor with a p-value greater than $0.05$ resulted in the parameter estimation in Table \ref{T3}. \\

\begin{table}[H]
\centerline{
\resizebox{0.48\textwidth}{!}{
\renewcommand{\arraystretch}{1.5}
\begin{tabular}{lcc}
\hline
\multicolumn{3}{c}{\bf{$\pi_h$ Coefficients}} \\
  \bf{Factor} &\bf{Coefficient}&\bf{P-Value}\\
  \hline
  $\log(\pi^m_{h})$&$1.08318$ &$< 2.2e-16$\\
  $\log(Pre)$&$0.42104$ &$0.02577$\\
  $\log(Class)$&$0.72842$ &$0.01093$\\
  \hline
\end{tabular}}}
\caption{Win Probability Coefficients} \label{T3}
\end{table}

The ~\cite{mcfad74} $R^2$ goodness of fit measure was used to compare the public's implied winning probabilities to the model's, where $R^2=1$ implies perfect predictive ability and $R^2=0$ means predictability is no better than random guessing. Using the last $30\%$ of the racing data, $R^2_{\pi_h}=0.218077$ and $R^2_{\pi^m_h}=0.214455$. We see the model has a small positive edge of $\Delta R^2=R^2_{\pi_h}-R^2_{\pi^m_h}=0.0036$ over the general public.

\subsection*{Estimating superfecta probabilities}

Below are the results of estimating the $\lambda^i$ parameters.

\begin{table}[H]
\centerline{
\resizebox{0.48\textwidth}{!}{
\renewcommand{\arraystretch}{1.5}
\begin{tabular}{lcc}
  \hline
  \multicolumn{3}{c}{\bf{Superfecta probability parameters}} \\
  \bf{Factor} &\bf{Coefficient}&\bf{P-Value}\\
  \hline
  $\lambda^1$&$0.600548$ &$< 2.2e-16$\\
  $\lambda^2$&$0.384509$ &$< 2.2e-16$\\
  $\lambda^3$&$0.26239$ &$7.767e-13$\\
  \hline
\end{tabular}}}
\caption{Superfecta probability parameters} \label{T4}
\end{table}

\subsection*{Estimating public's superfecta probabilities}

Below are the results of estimating the $\theta^i$ parameters.

\begin{table}[H]
\centerline{
\resizebox{0.48\textwidth}{!}{
\renewcommand{\arraystretch}{1.5}
\begin{tabular}{lcc}
  \hline
  \multicolumn{3}{c}{\bf{Superfecta probability parameters}} \\
  \bf{Factor} &\bf{Coefficient}&\bf{P-Value}\\
  \hline
  $\theta^1$&$1.2058$ &$< 2.2e-16$\\
  $\theta^2$&$0.8215$ &$< 2.2e-16$\\
  $\theta^3$&$0.5312$ &$< 2.2e-16$\\
  $\theta^4$&$0.4146$ &$< 2.2e-16$\\
  \hline
\end{tabular}}}
\caption{Superfecta probability parameters} \label{T4}
\end{table}

\end{document}